\documentclass[dvipdfmx]{article}
\usepackage{amsmath}
\usepackage{amssymb}
\usepackage{amsfonts}
\usepackage{amsthm}
\usepackage{blkarray}
\usepackage{arydshln}
\usepackage{multirow}
\usepackage[dvipdfmx]{graphicx}
\newcommand{\bigzerou}{%
\smash{\lower1.7ex\hbox{\bg 0}}}
\theoremstyle{definition}
\newtheorem{theorem}{Theorem}[section]

\newtheorem{lem}{Lemma}[section]

\newcommand{\ba}{\begin{eqnarray}}
\newcommand{\ea}{\end{eqnarray}}
\newcommand{\no}{\nonumber}

\newcommand{\mapright}[1]{%
\smash{\mathop{%
\hbox to 1.0cm{\rightarrowfill}}\limits^{#1}}}
\newcommand{\mapleft}[1]{%
\smash{\mathop{%
\hbox to 1.3cm{\leftarrowfill}}\limits^{#1}}}
\begin{document}
\title{
\begin{flushright}
  \begin{minipage}[b]{5em}
    \normalsize
    ${}$      \\
  \end{minipage}
\end{flushright}
{\bf Geometrical Proof of Generalized Mirror Transformation of Projective Hypersurfaces}}
\author{Masao Jinzenji\\
\\
\it Department of Mathematics, Okayama University\\
\it  Okayama, 700-8530, Japan\\
{\it e-mail address: pcj70e4e@okayama-u.ac.jp}}
\maketitle
\begin{abstract}
In this paper, we propose a geometrical proof of  the generalized mirror transformation of 
genus $0$ Gromov-Witten invariants of degree $k$ hypersurface in $CP^{N-1}$. 
\end{abstract}
\section{Introduction}
\subsection{Notation and Main Theorem}
Let $N,k$ be positive integers and $M_{N}^{k}$ be a degree $k$ hypersurface in $CP^{N-1}$. 
We denote by $\overline{M}_{0,n}(CP^{N-1},d)$ the moduli space of stable maps of degree 
$d$ from genus $0$ semi-stable curve with $n$ marked points to $CP^{N-1}$ \cite{K}.  
The genus $0$ $n$-pointed Gromov-Witten invariant of $M_{N}^{k}$ used in this paper is defined as follows  \cite{K, Jin2}, 
 \ba
\langle\prod_{j=1}^{n}{\cal O}_{h^{a_{j}}}\rangle_{0,d}=\int_{\overline{M}_{0,n}(CP^{N-1},d)}
c_{top}(\overline{\cal E}_{d}^{k})\wedge\biggl(\mathop{\bigwedge}_{j=1}^{n} ev_{j}^{*}(h^{a_{j}})\biggr),
\ea
where $h$ is the hyperplane class in $H^{*}(CP^{N-1},{\bf C})$ and $ev_{j}: \overline{M}_{0,2}(CP^{N-1},d)
\rightarrow CP^{N-1}\;\;(j=1,2,\cdots,n)$ is the evaluation map at the $j$-th marked point.
$\overline{\cal E}_{d}^{k}$ is the vector bundle on $\overline{M}_{0,2}(CP^{N-1},d)$ that impose the condition 
that image of the stable map is contained in the hypersurface.\footnote{This bundle is rigorously described by direct image sheaf of pull-back of ${\cal O}_{CP^{N-1}}(k)$ by 
evaluation map \cite{K, Jin2}, but we omit here this lengthy notation.} It is non-vanishing only if the following condition is satisfied,
\ba
\sum_{j=1}^{n}a_{j}=N-5+(N-k)d+n.
\ea

On the other hand, we also introduce the compactified moduli space of quasimaps from $CP^{1}$ with two marked points 
($0$ and $\infty$) to $CP^{N-1}$ of degree $d$, which we denote by $\widetilde{Mp}_{0,d}(N,d)$.
See \cite{Jin1} and \cite{S} for details of construction. We then define intersection number
$w({\cal O}_{h^{a}}{\cal O}_{h^{b}})_{0,d}$ of $\widetilde{Mp}_{0,d}(N,d)$, which is an analogue of $\langle{\cal O}_{h^{a}}{\cal O}_{h^{b}}\rangle_{0,d}$ of $\overline{M}_{0,2}(CP^{N-1},d)$,
as follows.
\ba
w({\cal O}_{h^{a}}{\cal O}_{h^{b}})_{0,d}=\int_{\widetilde{Mp}_{0,2}(N,d)}
c_{top}(\widetilde{\cal E}_{d}^{k})\wedge ev_{0}^{*}(h^a)\wedge ev_{\infty}^{*}(h^b).
\ea  
Here, $ev_{0}$ and $ev_{\infty}$ is the evaluation maps at $0$ and $\infty$ respectively and
$\widetilde{\cal E}_{d}^{k}$ is the vector bundle on $\widetilde{Mp}_{0,d}(N,d)$ which has the
same geometrical meaning as $\overline{\cal E}_{d}^{k}$. It is non-vanishing only if the following condition is satisfied.
\ba
a+b=N-3+(N-k)d.
\label{wsel}
\ea
In \cite{S}, Saito constructed explicit toric data of $\widetilde{Mp}_{0,d}(N,d)$ and showed that 
it is a compact toric orbifold. Moreover, he showed that its Chow ring is generated by $d+1$ 
divisor classes $H_{0},H_{1},\cdots,H_{d}$ that satisfy the following relations,
\ba
(H_{0})^{N}=0,\;(H_{i})^{N}(2H_{i}-H_{i-1}-H_{i+1})=0\;\;(i=1,2,\cdots,d-1),\;(H_{d})^{N}=0. 
\ea
In \cite{Jin1}, we showed that $w({\cal O}_{h^{a}}{\cal O}_{h^{b}})_{0,d}$
is written in terms of the Chow ring.
\ba
w({\cal O}_{h^{a}}{\cal O}_{h^{b}})_{0,d}=\int_{\widetilde{Mp}_{0,2}(N,d)}
(H_{0})^{a}\biggl(\frac{\prod_{j=1}^{d}e^{k}(H_{j-1}, H_{j})}{\prod_{j=1}^{d-1}(kH_{j})}\biggr)(H_{d})^{b},
\ea 
where $e^{k}(x,y)=\prod_{j=0}^{k}(jx+(k-j)y)$. In the above formula, we interpret that $kH_{j}$ in the denominator is canceled 
by the $kH_{j}$ in $e^{k}(H_{j-1}, H_{j})$. 
Therefore, we can compute $w({\cal O}_{h^{a}}{\cal O}_{h^{b}})_{0,d}$
explicitly.

Let $P_{g}$ be set of partitions of positive integer $g$:
\ba
P_{g}=\{\sigma_{g}=(g_{1},\cdots,g_{l(\sigma_{g})})\;|\;1\leq g_{1}\leq g_{2}\leq \cdots\leq g_{l(\sigma_{g})},\;\;
\sum_{j=1}^{l(\sigma_{g})}g_{j}=g\;\}.
\ea
For a partition $\sigma_{g}\in P_{g}$, we define multiplicity $\mbox{mul}(i,\sigma_{g})$ of $\sigma_{g}$ as follows.
\ba
\mbox{mul}(i,\sigma_{g})=(\mbox{number of subscript $j$ that satisfies $g_{j}=i$}).
\ea
We define combinatorial factor $S(\sigma_{g})$ as follows,
\ba
S(\sigma_{g})=\prod_{i=1}^{g}\frac{1}{(\mbox{mul}(i,\sigma_{g}))!}.
\ea

In this paper, we prove the following theorem that describes relation between intersection numbers $w({\cal O}_{h^{a}}{\cal O}_{h^{b}})_{0,d}$
and Gromov-Witten invariants $\langle\prod_{j=1}^{l}{\cal O}_{h^{a_{j}}}\rangle_{0,d}\;\;$. 
\begin{theorem}
\begin{eqnarray}
&&w({\cal O}_{h^{a}}{\cal O}_{h^{b}}
)_{0,d}-w({\cal O}_{h^{N-3-(k-N)d}}{\cal O}_{1}
)_{0,d} \no\\
&&=\langle{\cal O}_{h^{a}}{\cal O}_{h^{b}}\rangle_{0,d}+\no\\
&&+\sum_{g=1}^{d-1}\sum_{\sigma_{g}\in P_{g}}
S(\sigma_{g})
\langle{\cal O}_{h^{a}}{\cal O}_{h^{b}}
\prod_{i=1}^{l(\sigma_{g})}{\cal O}_{h^{1+(k-N)g_{i}}}\rangle_{0,d-g}\cdot\biggl(\prod_{i=1}^{l(\sigma_{g})}
\frac{w({\cal O}_{h^{N-3+(N-k)g_i}}{\cal O}_{1}
)_{0,g_i}}{k}\biggr),
\no\\
&& (a+b=N-3+(N-k)d).
\label{main}
\end{eqnarray}
\label{gt}
\end{theorem}
This theorem corresponds to the most concise version of ``generalized mirror transformation'' for K\"{a}hler sub-ring of the small quantum cohomology ring 
of $M_{N}^{k}$.

Generalized mirror transformation was first observed in \cite{Jin2} for low degrees and 
was generalized to arbitrary degree in \cite{Jin3} in the context of virtual structure constants \cite{Jin4}.  
Rigorous proof of generalized mirror transformation was given by Iritani \cite{Iri} by using Birkhoff factorization 
technique invented by Coates and Givental \cite{CG} and also by Guest \cite{Guest}.
But these results are based on three pointed Gromov-Witten invariants or $J$-function and process of explicit 
computation was quite complicated.

Later, we found that fundamental invariants to describe generalized mirror transformation are two pointed 
Gromov-Witten invariants. Therefore, we present here the version given in (\ref{main}) that is written 
in terms of $\langle{\cal O}_{h^{a}}{\cal O}_{h^{b}}\rangle_{0,d}$ and $w({\cal O}_{h^{a}}{\cal O}_{h^{b}})_{0,d}$.
As was mentioned before, we think that it is the most compact form of the generalized mirror transformation for 
small quantum cohomology ring. 
As for big quantum cohomology ring, generalized mirror transformation of projective hypersurfaces 
was conjectured in \cite{JS} and this version is much simpler for explicit computation of Gromov-Witten 
invariants. 

In this paper, we give a geometrical proof of the generalized mirror transformation for small quantum 
cohomology ring. The word ``geometrical'' means that we do not use localization technique.
Instead, we go back to our original motivation given in \cite{JN}.
That is to say, {\bf ``generalized mirror transformation is nothing but the process of removing contributions of
quasi maps that are not actual maps from $CP^{1}$ to $CP^{N-1}$''}.   
Our proof given in this paper clarifies geometrical meaning of the
generalized mirror transformation.  

\subsection{Usage of Theorem \ref{gt} and Historical Background of Quasimap }
We explain briefly usage of Theorem \ref{gt}. Let us first discuss the case of Fano hypersurface with $N-k\geq 2$.
Since $N-3-(k-N)d>N-2\;\;(d\geq 1)$, (\ref{main}) reduces to the following equality, 
\begin{eqnarray}
&&w({\cal O}_{h^{a}}{\cal O}_{h^{b}}
)_{0,d}=\langle{\cal O}_{h^{a}}{\cal O}_{h^{b}}\rangle_{0,d}.
\quad (a+b=N-3+(N-k)d).
\label{g22}
\end{eqnarray} 
This says that Gromov-Witten invariant is correctly evaluated by using the moduli space of quasimaps  $\widetilde{Mp}_{0,d}(N,d)$.
This fact was implied in \cite{Jin5} and follows from Theorem 9.1 in \cite{giv} proved by Givental. Explicit statement of the above equality 
was given in \cite{Jin4} in terms of the virtual structure constant:
\ba
\tilde{L}_{n}^{N,k,d}=\frac{d}{k}w({\cal O}_{h^{N-2-n}}{\cal O}_{h^{n-1+(N-k)d}})_{0,d}.
\label{vw}
\ea
(\ref{vw}) was proved for arbitrary $N$ and $k$ in \cite{Jin1}. 

If $N-k=1$, the above equality is slightly modified only in the $d=1$ case.
\ba
&&w({\cal O}_{h^{a}}{\cal O}_{h^{b}})-w({\cal O}_{h^{N-2}}{\cal O}_{1}
)_{0,1}=w({\cal O}_{h^{a}}{\cal O}_{h^{b}})-k\cdot k!
=\langle{\cal O}_{h^{a}}{\cal O}_{h^{b}}\rangle_{0,1}.
\quad (a+b=N-2).
\label{g2}
\ea 
It was also fundamentally proved in \cite{giv} and explicitly stated in \cite{CJ}. 

In the case of $N=k$ where the hypersurface is a Calabi-Yau hypersurface, 
we introduce the following generating function:
\ba
w({\cal O}_{h^{a}}{\cal O}_{h^{b}})_{0,d}(x):=kx+\sum_{d=1}^{\infty}w({\cal O}_{h^{a}}{\cal O}_{h^{b}})_{0,d}e^{dx}
\quad (a+b=N-3).
\ea
In \cite{Jin1}, we proved the following equality:
\ba
w({\cal O}_{h^{N-3}}{\cal O}_{1})_{0,d}(x)=kt(x),
\ea
where
\ba
&&t(x):=\bigl(x+\frac{w_{1}(x)}{w_{0}(x)}\bigr),\\
&&(w_{0}(x)=\sum_{d=0}^{\infty}\frac{(kd)!}{(d!)^{k}}e^{dx},\;\;w_{1}(x)=\sum_{d=1}^{\infty}\frac{(kd)!}{(d!)^{k}}(\sum_{i=1}^{d}\sum_{l=1}^{k-1}
\frac{l}{i(ki-l)})e^{dx}).
\ea
This $t(x)$ is nothing but the mirror map (in physics terminology, redefinition of the coupling constant of the Gauged Linear Sigma 
Model) used in the mirror computation of genus $0$ Gromov-Witten invariants of the Calabi-Yau hypersurface. Then Theorem \ref{gt}
is equivalent to the following equality:
\ba
w({\cal O}_{h^{a}}{\cal O}_{h^{b}})_{0,d}(x)=\langle{\cal O}_{h^{a}}{\cal O}_{h^{b}}\rangle_{0,d}(t(x)),
\label{rc}
\ea
where
\ba
\langle{\cal O}_{h^{a}}{\cal O}_{h^{b}}\rangle_{0,d}(t)
:=kt+\sum_{d=1}^{\infty}\langle{\cal O}_{h^{a}}{\cal O}_{h^{b}}\rangle_{0,d}e^{dt}.
\ea

Therefore by combining the results given in \cite{Jin1} and \cite{Jin4}, Theorem 1 gives a proof of the mirror theorem of genus $0$ 
Gromov-Witten invariants of the Calabi-Yau hypersurface. Of course, the mirror theorem in this case was already proved 
in \cite{giv} and \cite{lly}, but their treatment of the mirror transformation (\ref{rc}) is quite analytic or complicated. Therefore,
geometrical meaning of the mirror transformation from the point of view of quasimap is not clear in these works. Our proof 
of Theorem 1 provides a short and geometrically clear proof of the mirror theorem of the Calabi-Yau hypersurface. 

In the $N-k<0$ case where the hypersurface is general type,  Theorem \ref{gt} enables us to write down $\langle{\cal O}_{h^{a}}{\cal O}_{h^{b}}\rangle_{0,d}\; (a+b=N-3+(N-k)d)$ in terms of the virtual structure constant $\tilde{L}_{n}^{N,k,d^{\prime}}\;(d^{\prime}\leq d)$. The process is briefly 
given as follows. First, note that the equality:
\ba
d\langle{\cal O}_{h^{a}}{\cal O}_{h^{b}}\rangle_{0,d}=\langle{\cal O}_{h^{a}}{\cal O}_{h^{b}}{\cal O}_{h}\rangle_{0,d}.
\ea  
According to the reconstruction theorem of Kontsevich-Manin \cite{KM}, we can compute all the multi-point genus $0$
Gromov-Witten invariants $\langle\prod_{j=1}^{n}{\cal O}_{h^{a_j}}\rangle_{0,d^{\prime}}\;\; (d^{\prime}\leq d)$ from the 
initial data $\langle{\cal O}_{h^{a}}{\cal O}_{h^{b}}{\cal O}_{h}\rangle_{0,d^{\prime}} \;\;(d^{\prime}\leq d)$. Moreover,
in the $d=1$ case, Theorem \ref{gt} gives us,
\ba
\langle{\cal O}_{h^{N-2-n}}{\cal O}_{h^{n-1+(N-k)}}\rangle_{0,1}&=&w({\cal O}_{h^{N-2-n}}{\cal O}_{h^{n-1+(N-k)}}
)_{0,d}-w({\cal O}_{h^{N-3-(k-N)}}{\cal O}_{1})_{0,d}\no\\
&=&k(\tilde{L}_{n}^{N,k,1}-\tilde{L}_{1+k-N}^{N,k,1}).
\ea
Hence by induction of $d$, we can express $\langle{\cal O}_{h^{a}}{\cal O}_{h^{b}}\rangle_{0,d}\; (a+b=N-3+(N-k)d)$ in terms of the virtual structure constant $\tilde{L}_{n}^{N,k,d^{\prime}}\;(d^{\prime}\leq d)$ with the aid of the equality (\ref{vw}).
This procedure derives all the conjectures given in \cite{Jin2} and \cite{Jin3}, and completes the proof of the genus $0$ mirror theorem      
for general type hypersurfaces in our formulation.  The intersection number $w({\cal O}_{h^{N-3-(k-N)d}}{\cal O}_{1})_{0,d}$ also 
appears as the expansion coefficient of the mirror map of big quantum cohomology ring in the context of Iritani \cite{Iri}, but 
his proof is also quite analytic. Therefore, our proof clarifies geometrical meaning of the generalized mirror transformation also in this case. 

Lastly, we briefly review historical background of the theory of quasimap. The idea of quasimap was first introduced by Witten \cite{W}
and was also used by ourselves \cite{JN} in order to derive the leading term of the $(N-2)$-point genus $0$ correlation function of
Calabi-Yau hypersurface in $CP^{N-1}$. Then this idea was realized as Gauged Linear Sigma Model \cite{W, MP} in physics terminology 
and studied in the context of quantum field theory. This line of study was extended to the cerebrated work of Morrison and his collaborators
\cite{HKLM}. It seems that they succeeded in deriving hypergeometric series used in the mirror computation from 
the context of Gauged Linear Sigma Model,
but derivation of redefinition coupling constant was not done. Relation between the mirror transformation and the procedure of removing contributions from quasimaps that are not actual maps (in the context of \cite{LNS1,LNS2}, these are called freckled instantons) was qualitatively suggested in \cite{W, giv, LNS1,LNS2}, but these works lacks quantitative derivation of the mirror transformation.   
We think that key of quantitative derivation of the mirror transformation is the intersection number:
\ba
w({\cal O}_{h^{N-3-(k-N)d}}{\cal O}_{1})_{0,d},
\ea        
which does not vanish even with the trivial operator insertion ${\cal O}_{1}$. The reason why it does not vanish only comes from 
our construction of the moduli space $\widetilde{Mp}_{0,d}(N,d)$. This fact enables us to write down the following short proof of the 
generalized mirror transformation. Of course, $\widetilde{Mp}_{0,d}(N,d)$ is also independently constructed by Ciocan-Fontanine 
and Kim \cite{CK}. But their idea of construction is based on stability condition and quite different from our construction. Hence 
their derivation of the mirror transformation is quite different from ours. We hope that relation between the two different approaches 
of construction of the moduli space of quasi maps will be clarified in the future.   

As for mirror symmetry of general type hypersurface, Landau-Ginzburg model is considered as its mirror counterpart \cite{W2}. 
In \cite{GKR}, coincidence of hodge numbers of general hypersurface and corresponding Landau-Ginzbrug model is shown. Relation 
between intersection number of moduli space of quasi maps and Landau-Ginzburg model seems to be also important.   

This paper is organized as follows, In Section \ref{intromp}, we introduce the moduli space $\widetilde{Mp}_{0,2}(N,d)$, which plays the central role in the proof of Theorem \ref{gt}.
In Section 3, we prove Theorem \ref{gt}. In Appendix A, we explain the reason why the ``perturbation space'' introduced in Section 3 can be used to evaluate the contributions from the 
excess intersections to $w({\cal O}_{h^{a}}{\cal O}_{h^{b}})_{0,d}$.

\vspace{2cm}

{\bf Acknowledgment} 
We would like to thank Dr. Hayato Saito for discussions. 
Our research is partially supported by JSPS grant No. 22K03289.  

\section{The Moduli Space $\widetilde{Mp}_{0,2}(N,d)$}
\label{intromp}
Let ${\bf a}_{j},\;(j=0,1,\cdots,d)$ be vectors in ${\bf C}^{N}$ 
and let $\pi_{N}: {\bf C}^{N}\setminus \{{\bf 0}\} \rightarrow CP^{N-1}$ be a projection map. 
In this section, we define a degree $d$ quasimap $p$ from ${\bf C}^{2}$ to ${\bf C}^{N}$ as a map that consists of 
${\bf C}^{N}$vector-valued degree $d$ homogeneous polynomials in two coordinates $s,t$ of ${\bf C}^{2}$:    
\begin{eqnarray}
&&p:{\bf C}^{2}\rightarrow {\bf C}^{N}\no\\
&&p(s,t)={\bf a}_{0}s^{d}+{\bf a}_{1}s^{d-1}t+{\bf a}_{2}s^{d-2}t^{2}+\cdots+{\bf a}_{d}t^{d}.
\label{polyp1}
\end{eqnarray}
The parameter space of quasi maps is given by ${\bf C}^{N(d+1)}=\{ ({\bf a}_{0},{\bf a}_{1},\cdots,{\bf a}_{d}) \}$.
We denote by $Mp_{0,2}(N,d)$ the space obtained from dividing $\{({\bf a}_{0},\cdots,{\bf a}_{d})\in
{\bf C}^{N(d+1)}|\;{\bf a}_{0}\neq {\bf 0},\;{\bf a}_{d}\neq {\bf 0}\}$ by two ${\bf C}^{\times}$ actions induced from the following 
two ${\bf C}^{\times}$ actions on ${\bf C}^{2}$ via the map $p$ in (\ref{polyp1}). 
\begin{equation}
(s,t)\rightarrow ( \mu s,\mu t),\;\;(s,t)\rightarrow (s,\nu t). 
\label{torus1}
\end{equation}
With the above two torus actions, $Mp_{0,2}(N,d)$ can be regarded as the parameter space of degree $d$ quasi maps from $CP^{1}$ to $CP^{N-1}$ with
two marked points in $CP^{1}$: $0(=(1:0))$ and $\infty(=(0:1))$ . Set theoretically, it is given as follows:
\begin{eqnarray}
Mp_{0,2}(N,d) = \{ ( {\bf a}_{0},{\bf a}_{1},\cdots,{\bf a}_{d}) \in {\bf C}^{N(d+1)}\;|\;{\bf a}_{0},{\bf a}_{d}\neq {\bf 0}\}/({\bf C}^{\times})^2,
\end{eqnarray}
where the two ${\bf C}^{\times}$actions are given by, 
\begin{eqnarray}
&&( {\bf a}_{0},{\bf a}_{1},\cdots,{\bf a}_{d}) \rightarrow (\mu{\bf a}_{0},\mu{\bf a}_{1},\cdots,\mu{\bf a}_{d-1},\mu{\bf a}_{d})\no\\
&&( {\bf a}_{0},{\bf a}_{1},\cdots,{\bf a}_{d}) \rightarrow ({\bf a}_{0},\nu{\bf a}_{1},\cdots,\nu^{d-1}{\bf a}_{d-1},\nu^{d}{\bf a}_{d})
\label{two22}
\end{eqnarray} 
The condition ${\bf a}_{0}, {\bf a}_{d}\neq {\bf 0}$ assures that the images of $0$ and $\infty$ are well-defined in $CP^{N-1}$. 

At this stage, we have to note the difference between the moduli space of holomorphic maps from $CP^{1}$ to $CP^{N-1}$ and the moduli
space of quasi maps from $CP^{1}$ to $CP^{N-1}$. In short, the latter includes the points that are not actual maps from 
$CP^{1}$ to $CP^{N-1}$ but rational maps from $CP^{1}$ to $CP^{N-1}$. 
Let us consider a quasi map 
$\sum_{j=0}^{d}{\bf a}_js^{j}t^{d-j}$ which can be factorized as 
\begin{eqnarray}
\sum_{j=0}^{d}{\bf a}_js^{j}t^{d-j}=p_{d-d_1}(s,t)\cdot(\sum_{j=0}^{d_1}{\bf c}_js^{j}t^{d_1-j}),
\label{pfacin}
\end{eqnarray}
where $p_{d-d_1}(s,t)$ is a homogeneous polynomial of degree $d-d_1(>0)$ and $\sum_{j=0}^{d_1}{\bf c}_js^{j}t^{d_1-j})$ represents an actual 
holomorphic map of degree $d_{1}$ from $CP^{1}$ to $CP^{N-1}$. If we consider $\sum_{j=0}^{d}{\bf a}_js^{j}t^{d-j}$
as a map from $CP^{1}$ to $CP^{N-1}$, it should be regarded as a rational map whose images of the zero points of $p_{d-d_1}$ is 
undefined.Moreover, the closure of the image of this map is a rational curve of degree $d_1(<d)$ in $CP^{N-1}$. The reason why 
we include this kind of quasimap is that we can obtain simpler compactification of the moduli space than the moduli space of the stable maps 
$\overline{M}_{0,2}(CP^{N-1},d)$, the standard moduli space used to define the two-point Gromov-Witten invariants.  

Now, let us turn into the problem of compactification of $Mp_{0,2}(N,d)$. If $d=1$, $Mp_{0,2}(N,1)$ is given by,  
\begin{eqnarray}
Mp_{0,2}(N,1) = \{ ( {\bf a}_{0},{\bf a}_{1}) \in {\bf C}^{2N}\;|\;{\bf a}_{0},{\bf a}_{1}\neq {\bf 0}\}/({\bf C}^{\times})^2,
\end{eqnarray}
where $({\bf C}^{\times})^2$ action is given as follows.
\begin{eqnarray}
&&( {\bf a}_{0},{\bf a}_{1}) \rightarrow (\mu{\bf a}_{0},\mu{\bf a}_{1})\no\\
&&( {\bf a}_{0},{\bf a}_{1}) \rightarrow ({\bf a}_{0},\nu{\bf a}_{1}).
\label{two221}
\end{eqnarray} 
Therefore, $Mp_{0,2}(N,1)$ is nothing but $CP^{N-1}\times CP^{N-1}$ and is already compact. 
If $d\geq 2$, we have to use the two ${\bf C}^{\times}$ actions in (\ref{two22}) to turn ${\bf a}_0$ and ${\bf a}_d$ into the 
points in $CP^{N-1}$, $[{\bf a}_0]$ and $[{\bf a}_d]$. Therefore, we can easily see,
\begin{eqnarray}
Mp_{0,2}(N,d) = \{ ( [{\bf a}_{0}],{\bf a}_{1},\cdots,{\bf a}_{d-1},[{\bf a}_{d}]) \in CP^{N-1}\times {\bf C}^{N(d-1)}\times CP^{N-1}\;|\}/{\bf Z}_{d}.
\label{openmp}
\end{eqnarray} 
In (\ref{openmp}), the ${\bf Z}_d$ acts on ${\bf C}^{N(d-1)}$ as follows.
\begin{equation}
({\bf a}_{1},{\bf a}_{2}\cdots,{\bf a}_{d-1})\rightarrow ((\zeta_{d})^j{\bf a}_{1},(\zeta_{d})^{2j}{\bf a}_{2}\cdots,
(\zeta_{d})^{(d-1)j}{\bf a}_{d-1}), 
\end{equation} 
where $\zeta_d$ is the d-th primitive root of unity.  In this way, we can see that $Mp_{0,2}(N,d)$ is not compact if $d\geq 2$.
In order to compactify $Mp_{0,2}(N,d)$, we imitate the stable map compactification and add the following chains of quasi maps 
\begin{equation}
\cup_{j=1}^{l(\sigma_{d})}\bigl(\sum_{m_{j}=0}^{d_{j}-d_{j-1}} {\bf a}_{d_{j-1}+m_{j}}(s_{j})^{m_{j}}(t_{j})^{d_{j}-d_{j-1}-m_{j}}\bigr),
\;\;\bigl({\bf a}_{d_{j}}\neq {\bf 0},\;\;j=0,1,\cdots,l(\sigma_{d})\bigr),
\label{chain00}
\end{equation}
at the infinity locus of $Mp_{0,2}(N,d)$. In (\ref{chain00}), $d_{j}$'s are integers that satisfy,
\begin{equation}
1\leq d_1<d_2<\cdots<d_{l(\sigma_{d})}\leq d-1.
\end{equation}
We denote by $\widetilde{Mp}_{0,2}(N,d)$ the space 
obtained after this compactification. 
This $\widetilde{Mp}_{0,2}(N,d)$ is the moduli space we use in this paper.
It is explicitly constructed as a toric orbifold by introducing boundary divisor coordinates $u_1,u_2,\cdots 
u_{d-1}$ as follows. 
\begin{eqnarray}
&&\widetilde{Mp}_{0,2}(N,d) = \no\\
&&\{ ( {\bf a}_{0},{\bf a}_{1},\cdots,{\bf a}_{d},u_1,u_2,\cdots,u_{d-1}) \in {\bf C}^{N(d+1)+d-1}\;|\;{\bf a}_{0},
({\bf a}_{1},u_{1}),\cdots,({\bf a}_{d-1},u_{d-1}),{\bf a}_{d}\neq {\bf 0}\}/({\bf C}^{\times})^{d+1},\no\\
\end{eqnarray}
where the (d+1) ${\bf C}^{\times}$actions are given by,
\begin{eqnarray}
&&( {\bf a}_{0},{\bf a}_{1},\cdots,{\bf a}_{d},u_1,\cdots,u_{d-1}) \rightarrow (\mu_{0}{\bf a}_{0},\cdots,\mu_{0}^{-1}u_{1},
\cdots),\no\\
&& ( {\bf a}_{0},{\bf a}_{1},\cdots,{\bf a}_{d},u_1,\cdots,u_{d-1}) \rightarrow (\cdots,\mu_{1}{\bf a}_{1},\cdots,\mu_{1}^{2}
u_{1},\mu_{1}^{-1}u_{2},\cdots), \no\\
&& ( {\bf a}_{0},{\bf a}_{1},\cdots,{\bf a}_{d},u_1,\cdots,u_{d-1}) \rightarrow (\cdots,\mu_{i}{\bf a}_{i},\cdots,\mu_{i}^{-1}u_{i-1},\mu_{i}^{2}u_{i},
\mu_{i}^{-1}u_{i+1},\cdots),\;(i=2,\cdots,d-1),\no\\ 
&& ( {\bf a}_{0},{\bf a}_{1},\cdots,{\bf a}_{d},u_1,\cdots,u_{d-1}) \rightarrow (\cdots,\mu_{d-1}{\bf a}_{d-1},\cdots,\mu_{d-1}^{-1}u_{d-2},
\mu_{d-1}^{2}u_{d-1}
),\no\\ 
&& ( {\bf a}_{0},{\bf a}_{1},\cdots,{\bf a}_{d},u_1,\cdots,u_{d-1}) \rightarrow (\cdots,\mu_{d}{\bf a}_{d},\cdots,\mu_{d}^{-1}u_{d-1}). 
\label{action1}
\end{eqnarray}
In(\ref{action1}), "$\cdots$" in the r.h.s indicates that the ${\bf C}^{\times}$ actions are trivial.   
These torus actions are represented by a $(d+1)\times 2d$ weight matrix $W_{d}$:
\begin{eqnarray} 
W_{d}:=\bordermatrix{                    &{\bf a}_{0}&{\bf a}_{1}&{\bf a}_{2}&\cdots&{\bf a}_{d-3}&{\bf a}_{d-2}&{\bf a}_{d-1}&{\bf a}_{d}&u_{1}&
                                  u_{2}&u_{3}&\cdots&u_{d-2}&u_{d-1}\cr
                                  h_{0}&1&0&0&\cdots&0&0&0&0&-1&0&0&\cdots&0&0\cr 
                                  h_{1}&0&1&0&\cdots&0&0&0&0&2&-1&0&\cdots&0&0\cr  
                                  h_{2}&0&0&1&\ddots&0&\vdots&0&0&-1&2&-1&\ddots&0&0\cr   
                                  \vdots&\vdots&\vdots&\ddots&\ddots&\ddots&\vdots&\vdots&\vdots&\vdots&\ddots&\ddots&\ddots&\vdots&\vdots\cr 
                                  \vdots&\vdots&\vdots&0&\ddots&1&0&0&0&0&0&\ddots&\ddots&\ddots&0\cr  
                                  \vdots&\vdots&\vdots&\vdots&\ddots&0&1&0&0&0&0&\ddots&-1&2&-1\cr    
                                  h_{d-1}&0&0&0&\cdots&0&0&1&0&0&0&\cdots&0&-1&2\cr
                                  h_{d}&0&0&0&\cdots&0&0&0&1&0&0&\cdots&0&0&-1\cr}
\label{toric1}
\end{eqnarray}

Notice that the $A_{d-1}$ Cartan matrix appears in $W_{d}$.
If $u_{1},u_{2},\cdots,u_{d-1}\neq 0$, we can set all the $u_{i}$'s to 1 by using the $(d+1)$ torus actions. The remaining  two torus actions 
that leave them invariant are nothing but the ones given in (\ref{two22}). Therefore, the subspace given by the condition 
$u_{1},u_{2},\cdots,u_{d-1}\neq 0$ corresponds to $Mp_{0,2}(N,d)$.
If $u_{d_1}=0, u_{j}\neq 0\;\;(j\neq d_{1})$, we have to delete the $u_{d_1}$ column of matrix $W_{d}$. 
This operation turns the $A_{d-1}$ Cartan matrix into the $A_{d_{1}-1}\times A_{d-d_{1}-1}$ Cartan matrix and results in 
chains of two quasi maps:
\begin{equation}
(\sum_{j=0}^{d_{1}}{\bf a}_{j}s_{1}^{j}t_{1}^{d_{1}-j})\cup(\sum_{j=0}^{d-d_{1}}{\bf a}_{j+d_{1}}s_{2}^{j}t_{2}^{d-d_{1}-j}),\;\;
({\bf a}_{0},{\bf a}_{d_{1}},{\bf a}_{d}\neq {\bf 0}).
\end{equation}
Therefore, the corresponding boundary locus is given by $\displaystyle{Mp_{0,2}(N,d_1)\mathop{\times}_{CP^{N-1}}Mp_{0,2}(N,d-d_1)}$,
where $\displaystyle{\mathop{\times}_{CP^{N-1}}}$ is the fiber product with respect to the following projection maps:
\begin{eqnarray}
&&ev_{\infty}:Mp_{0,2}(N,d_1)\rightarrow CP^{N-1},\;\;ev_{\infty}({\bf a}_{0},\cdots,{\bf a}_{d_1})=[{\bf a}_{d_1}]\no\\ 
&&ev_{0}:Mp_{0,2}(N,d-d_1)\rightarrow CP^{N-1},\;\;ev_{0}({\bf a}_{d_1},\cdots,{\bf a}_{d})=[{\bf a}_{d_1}] 
\end{eqnarray}
In general, the subspace given by the condition 
\begin{equation}
u_{d_{i}}=0,\;(1\leq d_{1}<d_{2}<\cdots<d_{l(\sigma_{d})-1}\leq d-1), u_{j}\neq 0,\;\; (j\notin \{d_{1},d_{2},\cdots
,d_{l(\sigma_{d})-1}\}),
\end{equation} 
corresponds to 
chains of quasi maps labeled by ordered partition $\sigma_{d}=(d_{1}-d_{0},d_{2}-d_{1},d_{3}-d_{2},\cdots,d_{l(\sigma_{d})}-
d_{l(\sigma_{d})-1})$:
\begin{equation}
\cup_{j=1}^{l(\sigma_{d})}\bigl(\sum_{m_{j}=0}^{d_{j}-d_{j-1}} {\bf a}_{d_{j-1}+m_{j}}(s_{j})^{m_{j}}(t_{j})^{d_{j}-d_{j-1}-m_{j}}\bigr),
\;\;\bigl({\bf a}_{d_{j}}\neq {\bf 0},\;\;j=0,1,\cdots,l(\sigma_{d})\bigr),
\label{chain0}
\end{equation}
where we set $d_{0}=0,d_{l(\sigma_{d})}=d$. In this case, the corresponding boundary locus is,
\begin{equation}
Mp_{0,2}(N,d_1-d_0)\mathop{\times}_{CP^{N-1}}Mp_{0,2}(N,d_2-d_1)\mathop{\times}_{CP^{N-1}}\cdots\mathop{\times}_{CP^{N-1}}
Mp_{0,2}(N,d_{l(\sigma_d)}-d_{l(\sigma_d)-1}).
\end{equation} 
Since the lowest dimensional boundary:
\begin{equation}
Mp_{0,2}(N,1)\mathop{\times}_{CP^{N-1}}Mp_{0,2}(N,1)\mathop{\times}_{CP^{N-1}}\cdots\mathop{\times}_{CP^{N-1}}
Mp_{0,2}(N,1),
\end{equation} 
is identified with the compact space $(CP^{N-1})^{d+1}$, we can expect that $\widetilde{Mp}_{0,2}(N,d)$ is compact. As was mentioned before, 
the fact that $\widetilde{Mp}_{0,2}(N,d)$ is compact was proved by Saito \cite{S}.

\section{Proof of Theorem \ref{gt}}
As was shown in the previous section, $\widetilde{Mp}_{0,2}(N,d)$ has the following stratification.
\ba
&&\widetilde{Mp}_{0,2}(N,d)=\no\\
&&\coprod_{l=1}^{d}\coprod_{0=d_{0}<d_{1}<\cdots <d_{l-1}<d_{l}=d}
\biggl(Mp_{0,2}(N,d_1-d_0)\mathop{\times}_{CP^{N-1}}\cdots\mathop{\times}_{CP^{N-1}}
Mp_{0,2}(N,d_{l}-d_{l-1})\biggr).\no\\
\ea
Let us consider the stratum of highest dimension,
\ba
Mp_{0,2}(N,d)=\{({\bf a}_{0},{\bf a}_{1},\cdots,{\bf a}_{d})\;|\;{\bf a}_{i}\in{\bf C}^{N},{\bf a}_{0}\neq {\bf 0},\;{\bf a}_{d}\neq {\bf 0}\;\}
/({\bf C}^{\times})^{2}, 
\ea
where the two ${\bf C}^{\times}$ actions are given by,
\ba
&&({\bf a}_{0},{\bf a}_{1},\cdots,{\bf a}_{d-1},{\bf a}_{d})\rightarrow (\lambda^{d}{\bf a}_{0},\lambda^{d-1}{\bf a}_{1},\cdots,\lambda{\bf a}_{d-1},{\bf a}_{d}),\no\\
&&({\bf a}_{0},{\bf a}_{1},\cdots,{\bf a}_{d-1},{\bf a}_{d})\rightarrow ({\bf a}_{0},\nu{\bf a}_{1},\cdots,\nu^{d-1}{\bf a}_{d-1},\nu^{d}{\bf a}_{d}).
\label{equi}
\ea
$[({\bf a}_{0},{\bf a}_{1},\cdots,{\bf a}_{d-1},{\bf a}_{d})]$ represents a rational map $\varphi(s:t)=\pi_{N}(\sum_{j=0}^{d}{\bf a}_{j}s^{j}t^{d-j})$ from $CP^{1}$ to $CP^{N-1}$ modulo ${\bf C}^{\times}$ action on $CP^{1}$ that fixes $0=(0:1), \infty=(1:0)\in CP^{1}$:
\ba
(s:t)\rightarrow (s:\lambda t).
\label{p1}
\ea
Here, $\pi_{N}:{\bf C}^{N}\setminus\{{\bf 0}\}\rightarrow CP^{N-1}$ is the projective equivalence.
Note that $\sum_{j=0}^{d}{\bf a}_{j}s^{j}t^{d-j}$ is factorized into the following form up to ${\bf C}^{\times}$ 
multiplication.
\ba
\sum_{j=0}^{d}{\bf a}_{j}s^{j}t^{d-j}=\biggl(\prod_{j=1}^{l(\sigma_{g})}(\beta_{j}s-\alpha_{j}t)^{g_{j}}\biggr)
\cdot\biggl(\sum_{j=0}^{d-g}{\bf c}_{j}s^{j}t^{d-g-j}\biggr),
\label{factor}
\ea
where $\pi_{N}(\sum_{j=0}^{d-g}{\bf c}_{j}s^{j}t^{d-g-j})$ defines a well-defined map from $CP^{1}$ to $CP^{N-1}$.
 
Since ${\bf a}_{0},{\bf a}_{d}\neq {\bf 0}$, it follows that $z_{j}:=(\alpha_{j}:\beta_{j})$ never coincides with $0$
and $\infty$ in $CP^{1}$. These distinct points $(z_{j}\;(j=1,\cdots,l(\sigma_{g})))$ also represent points 
where $\varphi(s:t)$ is ill-defined. Since $[({\bf a}_{0},{\bf a}_{1},\cdots,{\bf a}_{d-1},{\bf a}_{d})]$ represents
$\varphi(s:t)$ modulo the ${\bf C}^{\times}$ action on $CP^{1}$, configuration of $z_{j}$'s should be considered modulo the
${\bf C}^{\times}$ action given in (\ref{p1}). 

We can easily see that the factorization in (\ref{factor}) is invariant under permutation of $z_{j}$'s with subscript 
$j$ that has the same value of $g_{j}$. With these considerations, we obtain the following 
decomposition of $Mp_{0,2}(N,d)$.
\ba
Mp_{0,2}(N,d)=\coprod_{g=0}^{d}\coprod_{\sigma_{g}\in P_{g}}
M_{0,2+l(\sigma_{g})}(CP^{N-1},d-g,\sigma_{g}).
\ea
In the above decomposition, $M_{0,2+l(\sigma_{g})}(CP^{N-1},d-g,\sigma_{g})$
is uncompactified moduli space of degree $d-g$ holomorphic map from $CP^{1}$ to $CP^{N-1}$ with $2+l(\sigma_{g})$
distinct marked points divided by $\prod_{i=1}^{g}\mbox{Sym}(\mbox{mul}(i,\sigma_{g}))$ action.
\ba
&&M_{0,2+l(\sigma_{g})}(CP^{N-1},d-g, \sigma_{g})\no\\
&&:=\{[\bigl(\pi_{N}(\sum_{j=0}^{d-g}{\bf c}_{j}s^{j}t^{d-g-j}),(0,\infty,z_{1},z_{2},\cdots,z_{l(\sigma_{g})})\bigr)]\}/\biggl(\prod_{i=1}^{g}\mbox{Sym}(\mbox{mul}(i,\sigma_{g}))\biggr).
\ea 
In the above definition, the tuple $\bigl(\pi_{N}(\sum_{j=0}^{d-g}{\bf c}_{j}s^{j}t^{d-g-j}),(0,\infty,z_{1},z_{2},\cdots,z_{l(\sigma_{g})})\bigr)$ is considered modulo the ${\bf C}^{\times}$ action
on $CP^{1}$ and the symmetric group $\mbox{Sym}(\mbox{mul}(i,\sigma_{g}))$ permutes $z_{j}$'s that satisfy $g_{j}=i$.  

Let ${\cal E}^{k}_{d}$ be a vector bundle on $Mp_{0,2}(N,d)$ that impose on $\varphi\in Mp_{0,2}(N,d)$ the condition 
that $\varphi(CP^{1})$ is contained in a generic degree $k$ hypersurface in $CP^{N-1}$. ${\cal E}^{k}_{d}$ is extended
to a rank $kd+1$ vector bundle $\widetilde{\cal E}^{k}_{d}$ on $\widetilde{Mp}_{0,2}(N,d)$.
Then we introduce the following intersection number of $\widetilde{Mp}_{0,2}(N,d)$.
\ba
w({\cal O}_{h^a}{\cal O}_{h^{b}})_{0,d}&=&\int_{\widetilde{Mp}_{0,2}(N,d)}c_{top}(\widetilde{\cal E}^{k}_{d})\wedge ev_{0}^{*}(h^{a})
\wedge ev_{\infty}^{*}(h^b).\no\\
&&(a+b=N-3+(N-k)d)
\ea
In the above definition, $ev_{0}$ and $ev_{\infty}$ are the evaluation maps at $(0:1)=0$ and $(1:0)=\infty$ respectively.
On the locus $M_{0,2+l(\sigma_{g})}(CP^{N-1},d-g, \sigma_{g})$, the condition that the section of $\widetilde{\cal E}^{k}_{d}$
induced from defining equation of the hypersurface vanishes, lowers the dimension only by $k(d-g)+1$. 
Hence if $k(d-g)+1+N-3+(N-k)d=N-2+Nd-kg\leq N-1+N(d-g)+l(\sigma_{g})-1\Longleftrightarrow (k-N)g+l(\sigma_{g})\geq 0$, 
this locus contributes to the intersection number. Moreover, if $(k-N)g+l(\sigma_{g})> 0$, we have excess intersection 
on the locus.
In order to evaluate contribution from this excess intersection, we introduce ``perturbation space'' defined as follows.
\ba
&&\widetilde{Mp}^{pert.}_{0,2}(N,d,\sigma_{g})\no\\
&&:=\left(M_{0,2+l(\sigma_{g})}(CP^{N-1},d-g)\biggl(\prod_{i=1}^{l(\sigma_{g})}(\mathop{\times}_{CP^{N-1}}
\widetilde{Mp}_{0,2}(N,g_{i}))\biggr)\right)/\biggl(\prod_{i=1}^{g}\mbox{Sym}(\mbox{mul}(i,\sigma_{g}))\biggr).
\ea 
where $M_{0,2+l(\sigma_{g})}(CP^{N-1},d-g)$ is the uncompactified moduli space of holomorphic maps 
from $CP^{1}$ with $2+l(\sigma_{g})$ marked points to $CP^{N-1}$,
\ba 
M_{0,2+l(\sigma_{g})}(CP^{N-1},d-g):=
\{[\bigl(\pi_{N}(\sum_{j=0}^{d-g}{\bf c}_{j}s^{j}t^{d-g-j}),(0,\infty,z_{1},z_{2},\cdots,z_{l(\sigma_{g})})\bigr)]\},
\ea
and the $i$-th fiber product is defined via the following diagrams,
\ba
M_{0,2+l(\sigma_{g})}(CP^{N-1},d-g)\stackrel{(ev_{i})_{i}}{\longrightarrow}(CP^{N-1})^{l(\sigma_{g})},
\label{dgm1}
\ea
and,
\ba
 \prod_{i}\widetilde{Mp}_{0,2}(N,g_{i})  \stackrel{\tiny{\prod_{i}ev_{0}}}{\longrightarrow} (CP^{N-1})^{l(\sigma_{g})}.
\ea
In the diagram (\ref{dgm1}), $ev_{i}$ is the evaluation map at $z_{i}$.
$\prod_{i=1}^{g}\mbox{Sym}(\mbox{mul}(i,\sigma_{g}))$ permutes $z_{j}$ together with $\mathop{\times}_{CP^{N-1}}\widetilde{Mp}_{0,2}(N,g_{j})$
in the same way as the definition of $M_{0,2+l(\sigma_{g})}(CP^{N-1},d-g, \sigma_{g})$.  

Let us explain geometrical meaning of $\widetilde{Mp}^{pert.}_{0,2}(N,d,\sigma_{g})$.
Roughly speaking, it generates ``infinitessimal deformation of quasimap in $M_{0,2+l(\sigma_{g})}(CP^{N-1},d-g,\sigma_{g})$'' that evaluates contribution of the excess intersection coming from
$M_{0,2+l(\sigma_{g})}(CP^{N-1},d-g,\sigma_{g})$, to $w({\cal O}_{h^a}{\cal O}_{h^{b}})_{0,d}$.\footnote{
The reason why this space can be used to evaluate contribution of the excess intersection will be explained in the proof of Lemma \ref{exs}.}         
In order to illustrate the idea, we introduce new homogeneous coordinates $(\tilde{s}:\tilde{t}):=(\beta_{i}s-\alpha_{i}t:t)$ of $CP^{1}$.
When $\sum_{j=0}^{d}{\bf a}_{j}s^{j}t^{d-j}$ is factorized in the form given in (\ref{factor}),
it is written in terms of the new homogeneous coordinates as follows:
\ba
\sum_{j=g_{i}}^{d}\tilde{\bf a}_{j}\tilde{s}^{j}\tilde{t}^{d-j}=\tilde{s}^{g_{i}}\biggl(\sum_{j=g_{i}}^{d}\tilde{\bf a}_{j}\tilde{s}^{j-g_{i}}
\tilde{t}^{d-j}\biggr),\quad (\tilde{\bf a}_{g_{i}},\;\tilde{\bf a}_{d}\neq{\bf 0}).
\label{defor1}
\ea
Therefore, deformation of this quasimap is given by,
\ba
\sum_{j=0}^{g_{i}}\tilde{\bf a}_{j}\tilde{s}^{j}\tilde{t}^{d-j}=\tilde{t}^{d-g_{i}}\biggl(\sum_{j=0}^{g_{i}}\tilde{\bf a}_{j}\tilde{s}^{j}\tilde{t}^{g_{i}-j}\biggr),
\label{defor2}
\ea  
which can be considered as a point in $\widetilde{Mp}_{0,2}(N,g_{i})$. Since $\tilde{\bf a}_{g_{i}}$ appears both
in (\ref{defor1}) and (\ref{defor2}), we use fiber product $\mathop{\times}_{CP^{N-1}}$. 
Note here that the quasimap $\sum_{j=0}^{g_{i}}\tilde{\bf a}_{j}\tilde{s}^{j}\tilde{t}^{g_{i}-j}$ should be
regarded as a quasimap from different $CP^{1}$ component attached to the original $CP^{1}$ at $z_{i}$, in order to describe ``infinitessimalness''
of the  deformation.
%Conventionally, we regard $(\tilde{s}:\tilde{t})=(1:0)$ as $0$ and  $(\tilde{s}:\tilde{t})=(0:1)$ as $\infty$
%on this component. 
Remaining construction will be given in Appendix A.   

Dimension of $\widetilde{Mp}^{pert.}_{0,2}(N,d,\sigma_{g})$ is given by,
\ba
&&N-1+N(d-g)-1+l(\sigma_{g})+\sum_{i=1}^{l(\sigma_{g})}(N-2+Ng_{i})-l(\sigma_{g})(N-1)\no\\
&&=N-2+Nd=\dim_{\bf C}(\widetilde{Mp}_{0,2}(N,d)).
\ea 
Hence this space can give us non-vanishing result. 
We can evaluate the contribution of the excess intersection coming from $M_{0,2+l(\sigma_{g})}(CP^{N-1},d-g,\sigma_{g})$
to the intersection number $w({\cal O}_{h^a}{\cal O}_{h^{b}})_{0,d}$ by extending the vector bundle $\widetilde{\cal E}^{k}_{d}$ to $\widetilde{Mp}^{pert.}_{0,2}(N,d,\sigma_{g})$. 
Then the contribution is given as follows,
\ba
&&\int_{\widetilde{Mp}^{pert.}_{0,2}(N,d,\sigma_{g})}c_{top}(\widetilde{\cal E}^{k}_{d})\wedge ev_{0}^{*}(h^{a})
\wedge ev_{\infty}^{*}(h^b) \no\\
&=&
\langle{\cal O}_{h^a}{\cal O}_{h^{b}}\prod_{i=1}^{l(\sigma_{g})}{\cal O}_{h^{1+(k-N)g_{i}}}\rangle_{0,d-g}
\biggl(\prod_{i=1}^{l(\sigma_{g})}\frac{w({\cal O}_{h^{N-3+(N-k)g_{i}}}{\cal O}_{1})_{0,g_{i}}}{k}\biggr)
\biggl(\prod_{i=1}^{g}\frac{1}{(\mbox{mul}(i,\sigma_{g}))!}\biggr),
\label{res}  
\ea
since  $w({\cal O}_{h^{N-3+(N-k)g_{i}}}{\cal O}_{1})_{0,g_{i}}$ (that contains insertion of ``$1$'' , or insertion of nothing) is non-vanishing. 
In deriving the above formula, we also used splitting axiom \cite{KM} of the fiber product:
\ba
&&\mathop{\times}_{CP^{N-1}}\; \mbox{under}\;\; M_{0,2+l(\sigma_{g})}(CP^{N-1},d-g)\stackrel{ev_{i}}{\longrightarrow}
CP^{N-1}\stackrel{ev_{0}}{\longleftarrow} \widetilde{Mp}_{0,2}(N,g_{i})\nonumber\\
&&\longleftrightarrow \frac{1}{k}\sum_{c_{i}=0}^{N-2}ev_{i}^{*}(h^{c_{i}})\wedge ev_{0}^{*}(h^{N-2-c_{i}}),
\ea
which follows from description of Poincar\'{e} dual of diagonal $\Delta$ in $CP^{N-1}\times  CP^{N-1}$ \;\;(precisely speaking $M_{N}^{k}\times M_{N}^{k}$) by cohomology elements.
Then we are naturally led to consider,
\ba
\sum_{c_{1}=0}^{N-2}\cdots \sum_{c_{l(\sigma_{g})}=0}^{N-2}\langle{\cal O}_{h^a}{\cal O}_{h^{b}}\prod_{i=1}^{l(\sigma_{g})}{\cal O}_{h^{c_{i}}}\rangle_{0,d-g}
\biggl(\prod_{i=1}^{l(\sigma_{g})}\frac{w({\cal O}_{h^{N-2-c_{i}}}{\cal O}_{1})_{0,g_{i}}}{k}\biggr),
\ea 
but we are forced to set $c_{i}=1+(k-N)g_{i}$ because we have the condition (\ref{wsel}).  

$\langle{\cal O}_{h^a}{\cal O}_{h^{b}}\prod_{i=1}^{l(\sigma_{g})}{\cal O}_{h^{1+(k-N)g_{i}}}\rangle_{0,d-g}$ in (\ref{res}) is the Gromov-Witten invariant of the hypersurface 
which is defined by,
\ba
&&\langle{\cal O}_{h^a}{\cal O}_{h^{b}}\prod_{i=1}^{l(\sigma_{g})}{\cal O}_{h^{1+(k-N)g_{i}}}\rangle_{0,d-g}\no\\
&&:=
\int_{M_{0,2+l(\sigma_g)}(CP^{N-1},d-g)}c_{top}({\cal E}^{k}_{d-g})\wedge ev_{0}^{*}(h^a)\wedge ev_{\infty}^{*}(h^b) \wedge\biggl(\mathop{\bigwedge}
_{i=1}^{l(\sigma_{g})}ev_{i}^{*}(h^{1+(k-N)g_{i}})\biggr).
\ea 
In the above formula, the moduli space $M_{0,2+l(\sigma_g)}(CP^{N-1},d-g)$ is not compactified. But note that boundary components added in compactification
by stable maps do not contribute to intersection numbers because they have positive codimensions under insertion of $c_{top}(\overline{\cal E}_{d-g}^{k})$. 
The last factor of the r.h.s. of  (\ref{res}) appears as the effect of dividing by the group $\prod_{i=1}^{g}\mbox{Sym}(\mbox{mul}(i,\sigma_{g}))$.

In the $d=g$ case, $\langle{\cal O}_{h^a}{\cal O}_{h^{b}}\prod_{i=1}^{l(\sigma_{g})}{\cal O}_{h^{1+(k-N)g_{i}}}\rangle_{0,0}$
vanishes if $l(\sigma_{g})$ is greater than $1$. 
If $l(\sigma_{g})=1$, we have,
\ba
 \langle{\cal O}_{h^a}{\cal O}_{h^{b}}{\cal O}_{h^{1+(k-N)d}}\rangle_{0,0}=k.
\ea
Of course, we can introduce perturbation spaces for the lower dimensional strata of $\widetilde{Mp}_{0,2}(N,d)$, but
they are irrelevant to the intersection number because they have positive codimensions.

In this way, we obtain the following equality:
\ba
&&w({\cal O}_{h^{a}}{\cal O}_{h^{b}})_{0,d}\no\\
&&=\langle{\cal O}_{h^{a}}{\cal O}_{h^{b}}\rangle_{0,d}\no\\
&&+\sum_{g=1}^{d-1}\sum_{\sigma_{g}\in P_{g}}\langle{\cal O}_{h^a}{\cal O}_{h^{b}}\prod_{i=1}^{l(\sigma_{g})}{\cal O}_{h^{1+(k-N)g_{i}}}\rangle_{0,d-g}
\biggl(\prod_{i=1}^{l(\sigma_{g})}\frac{w({\cal O}_{h^{N-3+(N-k)g_{i}}}{\cal O}_{1})_{0,g_{i}}}{k}\biggr)
\biggl(\prod_{i=1}^{g}\frac{1}{(\mbox{mul}(i,\sigma_{g}))!}\biggr)\no\\
&&+w({\cal O}_{h^{N-3+(N-k)d}}{\cal O}_{1})_{0,d}.
\ea
 This is nothing but the generalized mirror transformation given in (\ref{main}) !

 \vspace{3cm}

\appendix
\section{Some Comments on the Perturbation Space}
In this part, we add some comments on background idea of construction of the perturbation space:
\ba
&&\widetilde{Mp}^{pert.}_{0,2}(N,d,\sigma_{g})\no\\
&&=\left(M_{0,2+l(\sigma_{g})}(CP^{N-1},d-g)\biggl(\prod_{i=1}^{l(\sigma_{g})}(\mathop{\times}_{CP^{N-1}}
\widetilde{Mp}_{0,2}(N,g_{i}))\biggr)\right)/\biggl(\prod_{i=1}^{g}\mbox{Sym}(\mbox{mul}(i,\sigma_{g}))\biggr).
\ea 
As was suggested in (\ref{defor1}), the vector valued polynomial $\sum_{j=g_{i}}^{d}\tilde{\bf a}_{j}\tilde{s}^{j-g_{i}}
\tilde{t}^{d-j}$ defines a quasi map in $Mp_{0,2}(N,d-g_{i})$. If we regard $(\tilde{s}:\tilde{t})$ as the original homogeneous 
coordinates $(s:t)$ used in the construction of $\widetilde{Mp}_{0,2}(N,d)$, the quasi map $s^{g_{i}}\biggl(\sum_{j=g_{i}}^{d}\tilde{\bf a}_{j}s^{j-g_{i}}t^{d-j}\biggr)$ corresponds to the following boundary components of  $\widetilde{Mp}_{0,2}(N,d)$,
\ba
&&\coprod_{l=1}^{g_{i}}\coprod_{0=d_{0}<d_{1}<\cdots <d_{l-1}<d_{l}=g_{i}}
\biggl(Mp_{0,2}(N,d_1-d_0)\mathop{\times}_{CP^{N-1}}\cdots\mathop{\times}_{CP^{N-1}}
Mp_{0,2}(N,d_{l}-d_{l-1})\mathop{\times}_{CP^{N-1}}
Mp_{0,2}(N,d-g_{i})\biggr)\no\\
&&=\widetilde{Mp}_{0,2}(N,g_{i})\mathop{\times}_{CP^{N-1}}
Mp_{0,2}(N,d-g_{i}).
\ea
Of course, the part $\widetilde{Mp}_{0,2}(N,g_{i})\mathop{\times}_{CP^{N-1}}$ is not contained in the original 
$\widetilde{Mp}_{0,2}(N,d)$. But by applying $\widetilde{Mp}_{0,2}(N,g_{i})\mathop{\times}_{CP^{N-1}}$, we can create 
the perturbation without changing $\sum_{j=g_{i}}^{d}\tilde{\bf a}_{j}\tilde{s}^{j-g_{i}}\tilde{t}^{d-j}$ that produces the  corresponding boundary 
components of $\widetilde{Mp}_{0,2}(N,d)$ in case that we regard $z_{i}=(\alpha_{i}:\beta_{i})\;\;(\Longleftrightarrow 
(\tilde{s}:\tilde{t})=(0:1))$ as the $0=(0:1)$ in the original construction of $\widetilde{Mp}_{0,2}(N,d)$. 
Successive operation of $\widetilde{Mp}_{0,2}(N,g_{i})\mathop{\times}_{CP^{N-1}}\;\;(i=1,2,\cdots, l(\sigma_{g}))$ is the idea behind construction of the perturbation space. 

Next, we explain why the moduli space $\widetilde{Mp}_{0,2}(N,d)$ can evaluate the contribution from excess intersection. 
Let us take for example a subset of $Mp_{0,2}(N,d)$: 
\ba
U=\{\;[\tilde{\bf a}_{0}(s-t)^d](=[\tilde{\bf a}_{0}(\lambda s-t)^d]\;(\lambda\in {\bf C}^{\times}))\;|\;[\tilde{\bf a}_0]\in CP^{N-1}\}
\label{mexcess}
\ea 
, which corresponds to maximally degenerated locus. If we 
introduce new homogeneous coordinates $(\tilde{s}:\tilde{t})$ with $\tilde{s}=s-t$, the above quasimap is rewritten as 
$[\tilde{\bf a}_{0}\tilde{s}^d]$. Then the reason comes from the following lemma.
\begin{lem}
\label{exs}
The contribution from excess intersection coming from the subset $U$ in (\ref{mexcess}) to $w({\cal O}_{h^{a}}{\cal O}_{h^{b}})_{0,d}$
$(a+b=N-3+(N-k)d)$ equals to $w({\cal O}_{h^{N-3+(N-k)d}}{\cal O}_{1})_{0,d}$.  
\end{lem}
{\it proof)}  Since the image of the quasimap in $U$ is a point $[\tilde{\bf a}_{0}]\in CP^{N-1}$, the excess intersection
of  $w({\cal O}_{h^{a}}{\cal O}_{h^{b}})_{0,d}$ caused by $U$
is given by intersection of a codimension $a+b=N-3+(N-k)d$ hyperplane in $CP^{N-1}=\{[\tilde{\bf a}_{0}]\}$ and a degree $k$ hypersurface in the same $CP^{N-1}$. 

On the other hand, let us consider $w({\cal O}_{h^{N-3+(N-k)d}}{\cal O}_{1})_{0,d}$. It is defined by,  
\ba
w({\cal O}_{h^{N-3+(N-k)d}}{\cal O}_{1})_{0,d}=\int_{\widetilde{Mp}_{0,2}(N,d)}
c_{top}(\widetilde{\cal E}_{d}^{k})\wedge ev_{0}^{*}(h^{N-3+(N-k)d}).
\ea     
A point in $Mp_{0,2}(N,d)$ is represented by a quasimap:
\ba
{\bf a}_{0}s^d+{\bf a}_{1}s^{d-1}t+{\bf a}_{1}s^{d-2}t^2+\cdots+{\bf a}_{d}t^{d} \;\;({\bf a_{0}},{\bf a}_{1}\neq {\bf 0}),
\ea
where we can assume $||{\bf a}_{0}||=||{\bf a}_{d}||=1$ by using the two ${\bf C}^{\times}$ action in (\ref{equi}). 
But we have equivalence relation under the two ${\bf C}^{\times}$ action, and we can use another representative:
\ba 
{\bf a}_{0}s^d+(\lambda){\bf a}_{1}s^{d-1}t+(\lambda)^2{\bf a}_{1}s^{d-2}t^2+\cdots+(\lambda)^{d}{\bf a}_{d}t^{d} \;\;(||{\bf a_{0}}||=||{\bf a}_{1}||=1),
\ea
with fixed $\lambda\in {\bf C}^{\times}$. Then let us take the limit $\lambda\rightarrow 0$. This is possible because we have no operator insertion at $(s:t)=(0:1)$. Then the moduli space degenerates to $\{ [{\bf a}_{0}s^d ]\;|\;[{\bf a}_{0}]\in CP^{N-1}\;\}$
that equals to $U$ and the intersection represented by $w({\cal O}_{h^{N-3+(N-k)d}}{\cal O}_{1})_{0,d}$ reduces to  intersection of a codimension $N-3+(N-k)d$ hyperplane in $CP^{N-1}=\{[{\bf a}_{0}]\}$ and a degree $k$ hypersurface in the same $CP^{N-1}$.
Hence we obtain the assertion of the lemma.  $\Box$


\newpage

\begin{thebibliography}{99}
\bibitem{CG} T. Coates, A. B. Givental, \textit{Quantum Riemann?Roch, Lefschetz and Serre.} Ann. of Math. (2) 165 (1) (2007) 1--53. 
\bibitem{CJ} A.Collino, M.Jinzenji.  \textit{On the structure of the small quantum cohomology rings of projective hypersurfaces.} Comm. Math. Phys. 206 (1999), no. 1, 157--183.  
\bibitem{CK} I.Ciocan-Fontanine; B. Kim.
\textit{Wall-crossing in genus zero quasimap theory and mirror maps.} Algebr. Geom. 1 (2014), no. 4, 400--448. 
\bibitem{giv} A. B. Givental, \textit{Equivariant Gromov-Witten invariants.} Internat. Math. Res. Notices 1996, no. 13, 613--663. 
\bibitem{GKR} M. Gross, L. Katzarkov, H. Ruddat,  
\textit{Towards mirror symmetry for varieties of general type.} Adv. Math. 308 (2017), 208--275. 
\bibitem{Guest} M. A. Guest, \textit{From quantum cohomology to integrable systems.} Oxford Graduate Texts in Mathematics, 15. Oxford University Press, Oxford, 2008. ISBN: 978-0-19-856599-4.  
\bibitem{HKLM}H. Jockers, V. Kumar, J.-M.Lapan, D.-R.Morrison. \textit{Two-sphere partition functions and Gromov-Witten invariants. }
 Comm. Math. Phys. 325 (2014), no. 3, 1139--1170.
\bibitem{Iri} H. Iritani, \textit{Quantum D-modules and generalized mirror transformations.} Topology 47 (2008), no. 4, 225--276.
 \bibitem{Jin1} M. Jinzenji, \textit{Mirror Map as Generating Function of Intersection Numbers: Toric Manifolds with Two K\"{a}hler Forms}, arXiv:1006.0607, Comm. Math. Phys. 323(2013), no. 2, 747--811. 
\bibitem{Jin2} M. Jinzenji, \textit{On the quantum cohomology rings of general type projective hypersurfaces and generalized mirror transformation.} Internat. J. Modern Phys. A 15 (2000), no. 11, 1557--1595.
\bibitem{Jin3} M. Jinzenji, \textit{Coordinate change of Gauss-Manin system and generalized mirror transformation.} Internat. J. Modern Phys. A Ke20 (2005), no. 10, 2131--2156. 
\bibitem{Jin4}M. Jinzenji, \textit{Gauss-Manin system and the virtual structure constants.} Internat. J. Math. 13 (2002), no. 5, 445--477.
\bibitem{Jin5}M. Jinzenji, \textit{On quantum cohomology rings for hypersurfaces in $CP^{N-1}$.} J. Math. Phys. 38 (1997), no. 12, 6613--6638. 
\bibitem{JN} M. Jinzenji, M. Nagura, \textit{Mirror symmetry and an exact calculation of an $(N-2)$-point correlation function on a Calabi-Yau manifold embedded in $CP^{N-1}$.} Internat. J. Modern Phys. A 11 (1996), no. 7, 1217--1252.
\bibitem{JS} M. Jinzenji, M. Shimizu, \textit{Multi-Point Virtual Structure Constants and Mirror Computation of  $CP^{2}$-model.} Commun.Num.Theor Phys. 07 (2013) 411--468. 
\bibitem{K} M. Kontsevich. \textit{ Enumeration of rational curves via torus actions.  } {\it The moduli space of curves} (Texel Island, 1994), 335--368, Progr. Math., 129, Birkh\"{a}user Boston, Boston, MA, 1995.
\bibitem{KM} M. Kontsevich, Yu. I. Manin. \textit{Gromov-Witten classes, quantum cohomology, and enumerative geometry.} Comm. Math. Phys. 164 (1994), no. 3, 525--562. 
\bibitem{lly} B. H. Lian, K. Liu, S. -T. Yau. \textit{Mirror principle. I.} Asian J. Math. 1 (1997), no. 4, 729--763. 
\bibitem{LNS1} A. Losev, N. Nekrasov, S. Shatashvili. \textit{The Freckled Instantons.}  arXiv:hep-th/9908204.
\bibitem{LNS2} A. Losev, N. Nekrasov, S. Shatashvili. \textit{Freckled Instantons in Two and Four Dimensions.}  Class. Quant. Grav. 17: 1181--1187 (2000).   
\bibitem{MP}D. -R. Morrison, M. R. Plesser. \textit{Summing the Instantons: Quantum Cohomology and Mirror Symmetry in Toric Varieties.} Nucl. Phys. B 440: 279-354, (1995).
\bibitem{S} H. Saito, \textit{Chow Rings of $\widetilde{Mp}_{0,2}(N,d)$ and $\overline{M}_{0,2}({\bf P}^{N-1},d)$ and Gromov-Witten Invariants of Projective Hypersurfaces of Degree 1 and 2}, Internat. J. Math. 28(2017), no. 12, 1750090.
\bibitem{W} E. Witten. \textit{Phases of $N=2$ Theories In Two Dimensions.} Nucl. Phys. B403:159--222, (1993).
\bibitem{W2} E. Witten \textit{Dynamics of quantum field theory.} Quantum fields and strings: a course for mathematicians. Vol. 2. 
Edited by Pierre Deligne et al. American Mathematical Society, Providence, RI; Institute for Advanced Study (IAS), Princeton, NJ, 1119--1424 (1999). 
\end{thebibliography}
\end{document}